\documentclass[review]{elsarticle}
\usepackage[OT1]{fontenc}  	
\usepackage[spanish,french,english,USenglish]{babel} 			
\usepackage{amssymb}  			
\usepackage{amsmath}	
\usepackage{amsthm}			
\usepackage{esint}  				
\usepackage{bbm}   				
\usepackage{color,graphicx} 			
\usepackage[small]{caption}
\usepackage{mathpazo} 
\usepackage[colorlinks=true,linkcolor=webgreen,filecolor=webbrown,citecolor=webgreen]{hyperref}

\definecolor{webgreen}{rgb}{0,.5,0}
\definecolor{webbrown}{rgb}{.6,0,0}

\DeclareMathOperator*{\res}{res} 

\newcommand{\be}{\begin{equation}}
\newcommand{\ee}{\end{equation}} 

\newcommand{\seqnum}[1]{\href{http://oeis.org/#1}{\underline{#1}}}

\makeatletter
\newcommand{\specialnumber}[1]{
\def\tagform@##1{\maketag@@@{(\ignorespaces##1\unskip\@@italiccorr#1)}}}
\makeatother

\makeatletter
\def\ps@pprintTitle{%
\let\@oddhead\@empty
\let\@evenhead\@empty
\let\@oddfoot\@empty
\let\@evenfoot\@oddfoot
}\makeatother

\begin{document}

\begin{frontmatter}

\title{A Note on Some Recent Results for the Bernoulli Numbers \\ of the Second Kind\footnote{\phantom{-} \\[6mm]
\texttt{\underline{Note to the readers of the 3rd arXiv version:}
this version is a copy of the journal version of the article, which has been published in the 
Journal of Integer Sequences, vol.\ 20, no.\ 3, Article 17.3.8, pp.\ 1-7, 2017. 
URL: \url{https://cs.uwaterloo.ca/journals/JIS/VOL20/Blagouchine/blag5.html} \\
Artcile history: submitted 20 December 2016, accepted 26 January 2017, published 27 January 2017.\\ 
The layout of the present version and its page numbering differ from the journal version, but the content, the numbering of equations and the numbering 
of references are the same. 
For any further reference to the material published here, please, use the journal version of the paper, 
which you can always get for free from the journal site (Journal of Integer Sequences is an open access journal).}}}

\author{Iaroslav V.~Blagouchine\corref{cor1}} 
\ead{iaroslav.blagouchine@univ-tln.fr,
iaroslav.blagouchine@centrale-marseille.fr,
iaroslav.blagouchine@pdmi.ras.ru.}
\address{University of Toulon, France \\[1mm]
\& \\[1mm]
Steklov Institute of Mathematics at St.-Petersburg, Russia.}

\begin{abstract}
In a recent issue of the {\it Bulletin of the Korean Mathematical Society}, Qi and Zhang discovered  
an interesting integral representation for the Bernoulli numbers of the second kind (also known  
as {\it Gregory's coefficients}, {\it Cauchy numbers of the first kind},
and the {\it reciprocal logarithmic numbers}). 
The same representation also appears in many other sources,
either with no references to its author,
or with references to various modern researchers. 
In this short note, we show that this representation is a rediscovery of an old result obtained in the 19th century 
by Ernst Schr\"oder. 
We also demonstrate that 
the same integral representation may be readily derived by means of 
complex integration.  
Moreover, we discovered that the asymptotics of these numbers were also the subject of several rediscoveries,  
including very recent ones. 
In particular, the first-order asymptotics, which are usually (and erroneously) credited to Johan F.\ Steffensen, 
actually date back to the mid-19th century, and probably were known even earlier. 
\end{abstract}

\end{frontmatter}

\vskip.3in

\section{Rediscovery of Schr\"oder's integral formula} 
In a recent article in the {\it Bulletin of the Korean Mathematical So\-ci\-ety} \cite{qi_02},  
several results concerning the Bernoulli numbers of the second kind were presented.  
 
We recall that these numbers (OEIS \seqnum{A002206} and \seqnum{A002207}),  
which we denote below by $G_n$, are rational 
\begin{equation}\notag 
\begin{array}{llll} 
\displaystyle  
G_1\,=\,+\dfrac{1}{2} \,,\quad& G_2\,=\,-\dfrac{1}{12} \,,\quad& G_3\,=\,+\dfrac{1}{24} \,,\quad& G_4\,=\,-\dfrac{19}{720}\,, \\[5mm] 
G_5\,=\,+\dfrac{3}{160} \,,\quad& G_6\,=\,-\dfrac{863}{60480} \,,\quad& G_7\,=\,+\dfrac{275}{24192} \,,\quad& G_8\,=\,-\dfrac{33953}{3628800} \,,\quad \ldots 
\end{array} 
\end{equation} 
\vskip 0.07in  
\noindent and were introduced by the Scottish mathematician  
and astronomer James Gregory in 1670 in the context of area's interpolation formula.  
Subsequently, they were rediscovered by many famous mathematicians, including  
Gregorio Fontana, Lorenzo Mascheroni, Pierre-Simon Laplace, Augustin-Louis Cauchy, Jacques Binet, 
Ernst Schr\"oder, Oskar Schl\"omilch, Charles Hermite and many others. 
Because of numerous rediscoveries these numbers do not have a standard name, 
and in the literature they are also referred to as {\it Gregory's coefficients},  
{\it (reciprocal) logarithmic numbers}, {\it Bernoulli numbers of 
the second kind}, normalized {\it generalized Bernoulli numbers} $B_n^{(n-1)}$ and normalized {\it Cauchy numbers 
of the first kind} $C_{1,n}$. Usually, these numbers are defined either via their generating function 
\begin{eqnarray} 
\label{eq32} 
\frac{u}{\ln(1+u)} = 1+\sum 
_{n=1}^\infty G_n \, u^n ,\qquad|u|<1\,, 
\end{eqnarray} 
or explicitly 
\begin{eqnarray}\notag 
G_n\,=\,\frac{C_{1,n}}{n!}\,=\lim_{s\to n}\frac{-B_s^{(s-1)}}{\,(s-1)\,s!\,}= 
\,\frac{1}{n!}\! \int\limits_0^1\! x\,(x-1)\,(x-2)\cdots(x-n+1)\, dx\,,\qquad n\in\mathbb{N}\,. 
\end{eqnarray} 
It is well known that $G_n$ are alternating $\,G_n=(-1)^{n-1}|G_n|\,$ and decreasing in absolute value; 
they behave as $\,\big(n\ln^2 n\big)^{-1}\,$ at $n\to\infty$ and may be bounded from below and from above accordingly 
to formulas (55)--(56) from \cite{iaroslav_08}.  
For more information about these important numbers, see \cite[pp.\ 410--415]{iaroslav_08}, \cite[p.\ 379]{iaroslav_09}, 
and the literature given therein (nearly 50 references).

\begin{figure}[!t]    
\centering 
\includegraphics[width=0.83\textwidth]{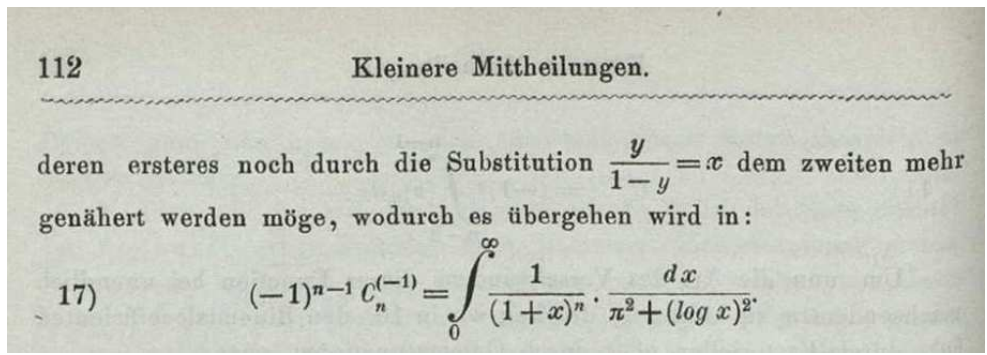} 
\caption{A fragment of p.\ 112 from Schr\"oder's paper \cite{schroder_01}. Schr\"oder's $C_n^{(-1)}$ are exactly our $G_n$.} 
\label{hkgcwuih} 
\end{figure} 
Now, the first main result of \cite[p.\ 987]{qi_02} is Theorem 1.\footnote{Our $G_n$ are exactly $b_n$ from \cite{qi_02} 
and $\frac{c_{n,1}^{(1)}}{n!}$ from \cite[Sect.\ 5]{chikhi_01}. Despite  
a venerable history, these numbers still lack a standard notation and various authors may use different notation for them.} It states:  
{\it the Bernoulli numbers of the second kind may be represented as follows} 
\begin{equation} 
G_n\,=\,(-1)^{n+1}\!\int\limits_1^\infty \!\!\frac{dt}{\big(\ln^2 (t-1)+\pi^2\big)\,t^n}\,,\qquad n\in\mathbb{N}\,. 
\end{equation} 
The same representation appears in a slightly different form\footnote{Put $t=1+u$.} 
\begin{equation}\label{jm23nd2d} 
G_n\,=\,(-1)^{n+1}\!\int\limits_0^\infty \!\!\frac{du}{\big(\ln^2 u +\pi^2\big)\,(u+1)^n}\,,\qquad n\in\mathbb{N}\,, 
\end{equation} 
in \cite[pp.\ 473--474]{coffey_02} and \cite[Sect.\ 5]{chikhi_01},
and is called {\it Knessl's representation}  
and {\it the Qi integral representation} respectively. 
Furthermore, various internet sources 
provide the same (or equivalent) formula, either with no references to its author or with references to different  
modern writers and/or their papers. 
However, the integral representation in question is not novel and is not due to Knessl nor to Qi and Zhang; in fact, 
this representation is a rediscovery of an old result. 
In a little-known paper of the German mathematician Ernst Schr\"oder \cite{schroder_01}, written in 1879, 
one may easily find exactly the same integral representation on p.\ 112;
see Fig.\ \ref{hkgcwuih}.  
Moreover, since this result is not difficult to obtain, it is possible that the same integral representation  
was obtained even earlier.

\section{Simple derivation of Schr\"oder's integral formula by means of the complex integration} 
Schr\"oder's integral formula \cite[p.\ 112]{schroder_01} may, of course, be derived in various ways. 
Below, we propose a simple derivation of this formula based on the method of contour integration. 
 
\begin{figure}[!t]    
\centering 
\includegraphics[width=0.35\textwidth]{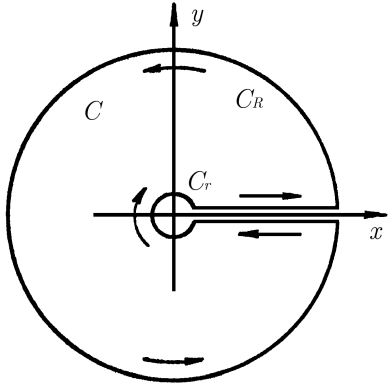} 
\caption{Integration contour $C$ ($r$ and $R$ are radii of the small and big circles respectively, where 
$r\ll1$ and $R\gg1$).} 
\label{kc30jfd} 
\end{figure}

If we set $u=-z-1$, then equality \eqref{eq32} may be written as 
\begin{eqnarray}\notag 
\frac{z+1}{\ln z - \pi i} = -1+\sum 
_{n=1}^\infty \big|G_n\big| \, (z+1)^n ,\qquad|z+1|<1\,. 
\end{eqnarray} 
Now considering the following line integral along a  
contour $C$ (see Fig.\ \ref{kc30jfd}),
where $n\in\mathbb{N}$, and then letting $R\to\infty\,$, $r\to0$, we have by the residue theorem 
\begin{eqnarray} 
&& \displaystyle  
\ointctrclockwise\limits_{C} \frac{\,dz\,}{\,(1+z)^n \, (\ln z - \pi i)\,}\,=\, 
\int\limits_{r}^R \! \ldots \, dz \, + \int\limits_{C_R} \!\ldots \, dz \,  
+ \int\limits_{R}^r \!\ldots \, dz \, + \int\limits_{C_r} \! \ldots \, dz \,  
\stackrel{\substack{R\to\infty \\ r\to0}}{=}   \notag\\[2mm] 
&& \displaystyle \quad 
= \, \int\limits_0^\infty \!\left\{\frac{1}{\ln x - \pi i} - \frac{1}{\ln x + \pi i} \right\}\cdot\frac{dx}{(1+x)^n} \, 
=\,2\pi i \! \int\limits_0^\infty \!\! \frac{1}{\,(1+x)^n} \cdot\frac{dx}{\,\ln^2 x +\pi^2\,}\,= \qquad\notag\\[3mm] 
&& \displaystyle  \quad  \notag 
=\,2\pi i \! \res\limits_{z=-1}\! \frac{1}{\,(1+z)^n \, (\ln z - \pi i)\,} \,=\,\frac{2\pi i}{n!}\cdot  
\left.\frac{d^n}{dz^n}\frac{z+1}{\,\ln z - \pi i\,}\right|_{z=-1} \!\!\! = \,2\pi i \, \big|G_n\big|\,, 
\end{eqnarray} 
since 
\begin{eqnarray} 
& \displaystyle  
\left| \,\int\limits_{C_R} \! \frac{\,dz\,}{\,(1+z)^n \, (\ln z - \pi i)\,} \, \right|  
\,=\, O\left(\!\frac{1}{\,R^{n-1}\ln R\,}\!\right)=o(1)\,, \qquad  & R\to\infty\,,\qquad n\geq1\,,\notag \\[5mm] 
& \displaystyle  
\left| \,\int\limits_{C_r} \! \frac{\,dz\,}{\,(1+z)^n \, (\ln z - \pi i)\,} \, \right|  
\,=\, O\left(\!\frac{r}{\,\ln r\,}\!\right)=o(1)\,, \qquad & 	r\to0\,,    \notag 
\end{eqnarray} 
and because at $z=-1$ the integrand of the contour integral has a pole of the $(n+1)$th order. This completes the proof. 
Note that above derivations are valid only for $n\geq1$, and so is Schr\"oder's integral formula, which 
may also be regarded as one of the generalizations of $G_n$ to the continuous values of $n.$

\section{Several remarks on the asymptotics for the Bernoulli numbers of the second kind} 
\begin{figure}[!t]    
\centering 
\includegraphics[width=0.5\textwidth]{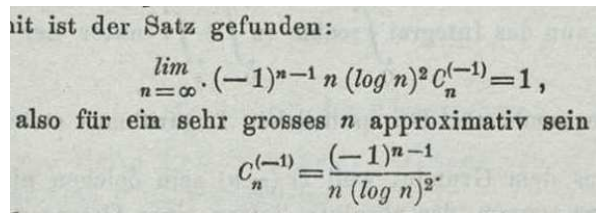} 
\caption{A fragment of p.\ 115 from Schr\"oder's paper \cite{schroder_01}.} 
\label{hg28g3dc} 
\vskip 0.25in  
\end{figure}

The first-order asymptotics $\,|G_n|\sim\big(n\ln^2 n\big)^{-1}\,$ at $n\to\infty$ are usually credited to  
Johan F.\ Steffensen \cite[pp.\ 2--4]{steffensen_01}, 
\cite[pp.\ 106--107]{steffensen_02}, \cite[p.\ 29]{norlund_01}, \cite[p.\ 14, Eq.\ (14)]{davis_02},  
\cite{nemes_01}, who found it in 1924.\footnote{The same first-order asymptotics also appear in \cite[p.\ 294]{comtet_01}, but without 
the source of the formula.}  
However, in our recent work \cite[p.\ 415]{iaroslav_08} we noted that exactly the same result appears in Schr\"oder's 
work written 45 years earlier, see Fig.\ \ref{hg28g3dc}, and the order of the magnitude of the closely related numbers  
is contained in a work of  
Jacques Binet dating back to 1839 \cite{binet_01}.\footnote{By the ``closely related numbers'' 
we mean the so-called {\it Cauchy numbers of the second kind} (OEIS 
\seqnum{A002657} and \seqnum{A002790}), and numbers $I'(k)$, see \cite[pp.\ 410--415, 428--429]{iaroslav_08}. The comment 
going just after Eq.\ (93) \cite[p.\ 429]{iaroslav_08} is based on the statements from \cite[pp.\ 231, 339]{binet_01}. 
The Cauchy numbers of the second kind $C_{2,n}$ and Gregory's coefficients $G_n$ are related to each other  
via the relationship $\,nC_{2,n-1}-C_{2,n}=n!\,|G_n|\,$, see \cite[p.\ 430]{iaroslav_08}.} In 1957 Davis \cite[p.\ 14, Eq.\ (14)]{davis_02} improved
this first-order  
approximation slightly by showing that $\,|G_n|\sim\Gamma(1+\xi) \big(n\ln^2 n + n\pi^2\big)^{-1}\,$ at $\,n\to\infty\,$ 
for some $\,\xi\in[0,1]\,$, without noticing that   
7 years earlier S.\ C.\ Van Veen had already obtained the complete asymptotics for them \cite[p.\ 336]{van_veen_01}, \cite[p.\ 29]{norlund_01}. 
Equivalent complete asymptotics were recently rediscovered in slightly different forms
by Charles Knessl \cite[p.\ 473]{coffey_02},  
and later by Gerg\H{o} Nemes  \cite{nemes_01}. 
An alternative demonstration of the same result was also presented by the author \cite[p.\ 414]{iaroslav_08}. 

\section{Acknowledgments}  
The author is grateful to Jacqueline Lorfanfant from the University of Strasbourg  
for sending a scanned version of \cite{schroder_01}.

\end{document}